\documentclass[12pt,leqno]{amsart}
\usepackage{latexsym}
\usepackage{amsthm}
\usepackage{amssymb}
\usepackage{amsmath}
\usepackage{color}

\newtheorem{thm}{Theorem}[section]
\newtheorem{lem}[thm]{Lemma}
\newtheorem{prop}[thm]{Proposition}
\newtheorem{cor}[thm]{Corollary}

\theoremstyle{remark}
    \newtheorem*{rem}{Remark}
\theoremstyle{definition}
    \newtheorem{defn}[thm]{Definition}

\numberwithin{equation}{section}
\def\Z {{\mathbb Z}}
\def\N {{\mathbb N}}
\def\Q {{\mathbb Q}}
\def\C {{\mathbb C}}
\def\T {{\tilde T}}
\def\e {{\tilde e}}
\def\H {{\mathcal H}}
\def\X {{\mathcal X}}
\def\G {\Gamma}

\title[the centre of the Hecke algebra]
{On the square root of the centre of the Hecke algebra of type $A$}
\author{Andrew Francis}
\address{School of Computing and Mathematics,\\ University of Western Sydney,\\ NSW 1797, Australia}
\author{Lenny Jones}
\address{Department of Mathematics,\\ Shippensburg University,\\ Pennsylvania, USA}

\begin{document}

\maketitle

\begin{abstract}
In this paper we investigate non-central elements of the Iwahori-Hecke algebra of the
symmetric group whose squares are central.

In particular, we describe a commutative subalgebra generated by certain non-central square roots of central elements, and the generic existence of a rank-three submodule of the Hecke algebra contained in the square root of the centre but not in the centre.
The generators for this commutative subalgebra include the longest word and elements related to trivial and sign characters of the Hecke algebra.
We find elegant expressions for the squares of such generators in terms of both the minimal basis of the centre and the elementary symmetric functions of Murphy elements.

\end{abstract}

\setcounter{section}{-1}

\section{Introduction}

In \cite{FranJones05OnBases}, we determined explicitly, for the centre of the Hecke algebra of type $A$, how to express each element of the $\Q[q,q^{-1}]$-norm basis in \cite{Jon90} as a linear combination of the elements of the $\Z[q,q^{-1}]$-minimal basis { \cite{Fmb}}. During our research for \cite{FranJones05OnBases}, we were led naturally to the square of the element of the Hecke algebra corresponding to the longest word in the symmetric group. This square, which arises in the definition of certain norms, is well-known to be central,
{and to have an expression as a product of Murphy elements which follows from its analogue in the braid group.  We derive directly a corresponding expression as a weighted sum of elementary symmetric functions of Murphy elements.}

{Motivated by recent conversations with John Murray, we
{have also investigated}
what other non-central elements of the Hecke algebra might have squares that are central.
(An analogous study by Murray for the centre of the symmetric
group algebra in the modular situation appears in \cite{Murray02}.) Results of a general nature arising from our investigation can be found in {Section \ref{sec:squareroot}},
while more detailed specific cases are handled in Section \ref{sec:other.sq.roots}.
The non-central elements we describe in Section \ref{sec:squareroot} are closely related to well-known central elements of the Hecke algebra, which are important in the study of its representation theory. Together with the longest word, these elements
span a submodule contained inside the square root of the centre,
and as with the element of the Hecke algebra corresponding to the longest word of the symmetric group, the squares of these elements can be
expressed as linear combinations of elementary symmetric functions
of Murphy operators {together with the identity} (Theorem \ref{thm:Ri^2=elsymfn}).}

\section{{Definitions} }\label{sec:Prelim}

Let $S=\{s_1,\dots,s_{n-1}\}$ be the standard generators of the symmetric group $S_n$, meaning that $s_i=(i\ i+1)$ as a permutation.
We say an expression $w = s_{i_1}s_{i_2}\dots s_{i_k}\in S_n$, $s_{i_j}\in  S$, is \emph{reduced} if $k$ is minimal. In this case, we say the length of $w$, denoted
by $l(w)$, is $k$.
A \emph{partition} of $n$ is a composition whose components are weakly decreasing from left to right. If $\lambda$ is a partition of $n$ we write $\lambda\vdash n$. The conjugacy classes in $S_n$ are indexed by partitions of $n$.

There are two fairly standard and closely related sets of generators and relations with which to define {the Hecke algebra $\H:=\H_n$ of $S_n$. Let $q$ be an indeterminate and let {$R=\Z[q^{1/2},q^{-1/2}]$}. Then $\H$ is the associative $R$-algebra generated by the set $\{T_s\mid s\in S\}$ with identity $T_1$ and subject to the relations:}

\begin{align}
  T_s^2&=qT_1+(q-1)T_s,     &&\text{ for }s\in S\label{eq:rel.q-order}\\
  T_{s_i}T_{s_{i+1}}T_{s_i}&=T_{s_{i+1}}T_{s_i}T_{s_{i+1}},&&\text{ for }1\le i\le n-2\label{eq:rel.q-braid}\\
  T_{s_i}T_{s_j}&=T_{s_j}T_{s_i}&&\text{ for }|i-j|\ge 2.\label{eq:rel.q-braid.commute}
\end{align}
The algebra $\H$ is a free $R$-module with basis $\{T_w:=T_{s_{i_1}}\dots T_{s_{i_r}}\mid w\in S_n\}$ where $w=s_{i_1}\dots s_{i_r}$ is a reduced expression for $w$ in $S_n$. {We present} most of the results in this paper in terms of this set of generators and relations, and with respect to this basis of $\H$ as a module. {However, some results (Section \ref{sec:sq.longword})} are more elegantly expressed in terms of a normalised generating set $\{\T_s\mid s\in S\}$, defined by {setting $\T_s:=q^{-1/2}T_s$}.  If we set $\xi=q^{1/2}-q^{-1/2}$ then the new generators satisfy the braid relations \eqref{eq:rel.q-braid} {and \eqref{eq:rel.q-braid.commute}}, but the order relation \eqref{eq:rel.q-order} becomes
\begin{equation}
\T_s^2=\T_1+\xi\T_s.\label{eq:rel.xi-order}
\end{equation}
We indicate corresponding elements and substructures of this normalised version of the Hecke algebra by placing a tilde above them. For example, we denote the normalised Hecke algebra itself as $\tilde{\H}$. Then $\tilde\H$ is a free $\Z[\xi]$-module with basis $\{\T_w:=\T_{s_{i_1}}\dots \T_{s_{i_r}}\}$ where $w=s_{i_1}\dots s_{i_r}$ is a reduced expression for $w$ in $S_n$.

The symmetric group $S_n$ has a symmetry best described as the
graph automorphism induced by reflecting the Dynkin diagram {about}
its midpoint.  If $\rho_n:S_n\to S_n$ is the group automorphism
defined by this symmetry then $\rho_n(s_i)=s_{n-i}$ for $1\le i\le
n-1$.  This automorphism is naturally extended to an algebra
{automorphism of $\H$}, also denoted $\rho_n$.

The \emph{centre} {$Z$ of $\H$} is defined to be the set of elements $c\in\H$ such that $ch=hc$ for all $h\in\H$.  {Throughout this paper we make frequent use of the `minimal basis' for $Z$ given in (\cite{GR97,Fmb}). This minimal basis} $\{\Gamma_\lambda\mid\lambda\vdash n\}$ {is an $R$}-basis for $Z$ and can be characterized by the following two properties:
\begin{itemize}
  \item $\Gamma_\lambda|_{q=1}=\sum_{w\in C_\lambda}T_w$, {where $C_\lambda$ is the conjugacy class of cycle type $\lambda$ in $S_n$}, and
  \item $\Gamma_\lambda-\sum_{w\in C_\lambda}T_w$ contains no minimal length elements from any conjugacy class.
\end{itemize}
The corresponding {class elements in $\tilde\H$ are}
$\{\tilde\Gamma_\lambda\mid\lambda\vdash n\}$ with the characterizing properties $\tilde\Gamma_\lambda|_{\xi=0}=\sum_{w\in C_\lambda}\T_w$, and $\tilde\Gamma_\lambda-\sum_{w\in C_\lambda}\T_w$ contains no minimal length elements from any conjugacy class.
We have $\tilde\Gamma_\lambda=q^{-l_\lambda/2}\Gamma_\lambda$, where $l_\lambda$ is the length of the shortest elements in $C_\lambda$.

Write $w_n$ for the longest word in $S_n$.  Define {$\ell :=l(w_n)=\sum_{i=1}^{n-1}i$}, the length of the longest word.

Define the \emph{Poincar\'e polynomial} to be
\[p(q):= \sum_{w\in S_n}q^{l(w)}.\]

\section{ {The square of the element in $\tilde\H$ corresponding to $w_n$} }\label{sec:sq.longword}

In this section we use the normalized version of the {Hecke algebra $\tilde\H$}
to give an explicit and surprisingly simple expression for $\T_{w_n}^2$.

We begin with some basic facts about Murphy elements.

\begin{defn} The \emph{Murphy operators} or \emph{Murphy elements} {$L_{n,i}$ of {$\H$} are defined by setting $L_{n,1}=0$ and
\begin{align*}
L_{n,i}&=T_{(i-1\ i)}+q^{-1}T_{(i-2\ i)}+\dots+q^{-(i-2)}T_{(1\ i)}\\
&=T_{s_{i-1}}+q^{-1}T_{s_{i-2}s_{i-1}s_{i-2}}+\dots+q^{-(i-2)}T_{s_1\dots s_{i-1}\dots s_1}
\end{align*}
for $i=2,3,\dots, n$.}
\end{defn}

{By a \emph{symmetric function in a set of commuting variables}, we mean a polynomial in those variables that is unchanged by any permutation of the variables.}

\begin{thm}[Murphy \cite{Mur83} or Dipper-James \cite{DJ87}]\label{thm:murphys.central} The Murphy elements commute.  The symmetric functions in the Murphy elements are in $Z$.
\end{thm}

The following was shown by Dipper and James \cite{DJ87} when {the Hecke algebra} is defined over a suitable field, and conjectured by them to be true over $R$.
Let {$\H_{\Q(q^{1/2})}$ be the Hecke algebra taken over $\Q(q^{1/2})$, and let $Z_{\Q(q^{1/2})}$ be the centre of $\H_{\Q(q^{1/2})}$.  Elements of $\H$ (such as the Murphy elements) are naturally embedded in $\H_{\Q(q^{1/2})}$.}

\begin{thm}[Dipper-James \cite{DJ87}]
The set of symmetric functions in the Murphy elements in {$\H_{\Q(q^{1/2})}$ is precisely $Z_{\Q(q^{1/2})}$.}
\end{thm}

The Murphy elements can be similarly defined {in $\tilde\H$}
by setting  $\tilde L_{n,1}:=0$ and
$\tilde L_{n,i}:=\T_{(i-1\ i)}+\T_{(i-2\ i)}+\dots+\T_{(1\ i)}$.  Of course, the statement of Theorem \ref{thm:murphys.central} readily translates to the normalised context.

The $i${-}th \emph{elementary symmetric function} in $n$ commuting variables is the sum of all monomials of length $i$ in the variables whose exponents are at most 1. For example, the 2nd elementary symmetric function in {the four} variables {$x_1,x_2,x_3,x_4$} is $x_1x_2+x_1x_3+x_1x_4+x_2x_3+x_2x_4+x_3x_4$.
Define {$\e_{n,0}:=\T_1$}, and for $i=1,2,\dots , n-1$, let {$\e_{n,i}$} denote the $i$-th elementary symmetric function in the Murphy elements $\tilde L_{n,1},\dots,\tilde L_{n,n}$.

\begin{prop}[\cite{FG:DJconj}]\label{thm:bauble}
For $i=0,\dots,n-1$,
\[ \e_{n,i} =\sum_{l_\lambda=i}\tilde\Gamma_\lambda.\]
Moreover, {we have the corresponding} $e_{n,i}=\sum_{l_\lambda=i}\Gamma_\lambda$.
\end{prop}

\begin{cor} The set {$\{\sum_{l_\lambda=i}\Gamma_\lambda\mid 1\le i\le n-1\}$} generates {$Z_{\Q(q^{1/2})}$} over {$\Q(q^{1/2})$}.
\end{cor}

Set
\begin{align*}  {\tilde M}_{n,i}: & =\rho_n(\tilde L_{n,i}) \\
                           &=\T_{(n\ n-i+1)}+\dots+\T_{(n-i+2\ n-i+1)},
\end{align*}
for $1\le i\le n$.  We call the {$\tilde M_{n,i}$} the \emph{dual Murphy elements}. These {elements} also appear in {\cite[p.26]{OgievPyatov99}}.

For example, in $S_4$

\begin{align*}
\tilde M_{4,1}&= {0}  &&\\
\tilde M_{4,2}&=\T_{(4\ 3)}&&=\T_{s_3}\\
\tilde M_{4,3}&=\T_{(4\ 2)}+\T_{(3\ 2)}&&=\T_{s_2s_3s_2}+\T_{s_2}\\
\tilde M_{4,4}&=\T_{(4\ 1)}+\T_{(3\ 1)}+\T_{(2\ 1)}&&=\T_{s_1s_2s_3s_2s_1}+\T_{s_1s_2s_1}+\T_{s_1}.
\end{align*}

These {elements} are `dual' because if you add `down the columns' instead of along the rows you get the standard Murphy elements.  {That is, adding the first terms of each element one has $\tilde L_{4,4}$, adding the second terms gives $\tilde L_{4,3}$, and {adding} the third terms gives $\tilde L_{4,2}$}.

If you like, the $i${-}th Murphy element `hangs down' from $i$ to 1, while the dual version `grows up' from $n-i+1$ to $n$.

\begin{lem}\label{lem:dualMeltsbasicprops}\text{}
\begin{enumerate}
\item {$\sum_{i=1}^n\tilde M_{n,i}=\sum_{i=1}^n\tilde L_{n,i}$}
\item\label{lem:barLi:i-1} {$\tilde M_{i,i}=\tilde M_{i-1,i-1}+\T_{(1\ i)}$ for $i\ge 2$}
\item\label{lem:T1-nTn-1} {$\T_{s_1\dots s_n}\T_{s_n\dots s_1}= { \T_1 } +\xi\tilde M_{n+1,n+1}$}
\item\label{lem:en.and.Ln+1.to.en+1} For {$i=0,1,\dots ,n-2$}, {$\tilde L_{n+1,n+1} \tilde e_{n,i}+\tilde e_{n,i+1}=\tilde e_{n+1,i+1}$.}
For $i=n-1$ we have {$\tilde L_{n+1,n+1}\tilde e_{n,n-1}=\tilde e_{n+1,n}$}.
\item\label{lem:rho.fixes.elsymfns} {$\rho_{n}(\tilde e_{n,i})=\tilde e_{n,i}$} {for $i=0,1,\dots , n-1$.}
\end{enumerate}
\end{lem}

\begin{proof}\text{}
\begin{enumerate}
\item {The sums on each side are equal to $\displaystyle\sum_{1\le i<j\le n}\T_{(i\ j)}$}.
\item This follows by definition of {$\tilde M_{i,j}$}.
\item This involves a straightforward induction on $n$ using \eqref{lem:barLi:i-1}.
\item These are elementary properties of symmetric functions in commuting variables.
\item This also follows since {$\tilde e_{n,i}$} is symmetric not only in the Murphy elements but also in the $s_i$ under $\rho_n$.
\end{enumerate}
\end{proof}

The following result {(Eq. \eqref{eq:sq.Ogiev-Pyatov})} {is stated in \cite[\S3]{Ram97} and \cite[\S 6]{OgievPyatov99}} using an alternative definition of Murphy elements {in $\tilde\H$, defined by setting} $\mathcal L_{n,1}=1$ and $\mathcal L_{n,i+1}=\T_{s_i}\mathcal L_{n,i}\T_{s_i}$.  It can be proved by induction on $n$, but is also a direct consequence of looking at the Hecke algebra as a quotient of the Braid group algebra.  In the braid group, the square of the longest word is the `full twist', and the Murphy element $\mathcal L_{n,i}$ is the braid whose $i$-th string goes behind strings $i-1,\dots,1$ and returns in front of strings $1,\dots,i-1$ to the $i$-th position.  All other strings remain as in the identity braid.  Equation \eqref{eq:sq.Ogiev-Pyatov} follows immediately.
\begin{equation}\label{eq:sq.Ogiev-Pyatov}{\T_{w_n}^2=\prod_{i=1}^{n} { {\mathcal L}_{n,i} }. } \end{equation}
We {present in Theorem \ref{thm:Twn^2=sumxi_ie_i}} an alternative expression for $\T_{w_n}^2$ which is a direct consequence of \eqref{eq:sq.Ogiev-Pyatov} together with the relation {$\mathcal L_{n,i}=\xi\tilde L_{n,i}+\T_1$} (the proof of this relation is elementary).  We provide a direct independent proof here.
\begin{thm}\label{thm:Twn^2=sumxi_ie_i} {$\displaystyle\T_{w_n}^2=\sum_{i=0}^{n-1}\xi^i\tilde e_{n,i}$.}
\end{thm}

\begin{proof}
{The proof is by induction on $n$.}  If $n=2$ then $\T_{w_2}^2=\T_{s_1}^2= { \T_1 } +\xi\T_{s_1}=\tilde e_{2,0}+\xi \tilde e_{2,1}$ since $\tilde e_{2,1}=\T_{s_1}$.

Suppose the statement is true for $n$.

Now $w_{n+1}=s_1\dots s_n w_n=w_ns_n\dots s_1$.  Also, it is clear that $w_ns_n\dots s_1=s_n\dots s_1\rho_{n+1}(w_n)$.  So, {using the fact that $\rho_{n+1}$ is an algebra automorphism,}
\begin{align*}
\T_{w_{n+1}}^2  &=\T_{s_1\dots s_n}\T_{w_n}^2\T_{s_n\dots s_1}\\
                &=\T_{s_1\dots s_n}\T_{s_n\dots s_1}\rho_{n+1}\left(\T_{w_n}^2\right)\\
                &{ =\rho_{n+1}\left(\mathcal L_{n+1,n+1}\right)\rho_{n+1} \left(\sum_{i=0}^{n-1}\xi^i\tilde e_{n,i}\right)\quad \text{by induction} } \\
                &=\rho_{n+1}\left(\left(\T_1 +\xi \tilde L_{n+1,n+1}\right)\left(\sum_{i=0}^{n-1}\xi^i\tilde e_{n,i}\right)\right)\\
                &=\rho_{n+1}\bigl( \T_1 +\xi( \tilde L_{n+1,n+1}\tilde e_{n,0}+\tilde e_{n,1})+\xi^2( \tilde L_{n+1,n+1}\tilde e_{n,1}+\tilde e_{n,2})+\dots\\
                &\qquad+\xi^{n-1}( \tilde L_{n+1,n+1}\tilde e_{n,n-2}+\tilde e_{n,n-1})+\xi^n\tilde L_{n+1,n+1}\tilde e_{n,n-1}\bigr)\\
                &=\rho_{n+1}\bigl( \T_1 +\xi \tilde e_{n+1,1}+\xi^2\tilde e_{n+1,2}+\dots+\xi^n\tilde e_{n+1,n}\bigr)\\
                &= \T_1 +\xi \tilde e_{n+1,1}+\xi^2\tilde e_{n+1,2}+\dots+\xi^n\tilde e_{n+1,n}
\end{align*}
as desired (the last equality due to Lemma \ref{lem:dualMeltsbasicprops}(\ref{lem:rho.fixes.elsymfns}), the second {to} last due to \ref{lem:dualMeltsbasicprops}(\ref{lem:en.and.Ln+1.to.en+1})).
\end{proof}

{
\begin{cor}\label{cor:sq.wn.gamma}
  In terms of class elements,
  \[\T_{w_n}^2=\sum_{\lambda\vdash n}\xi^{l_\lambda}\tilde\Gamma_\lambda.\]
\end{cor}
}
{
\begin{proof}
  This is an immediate consequence of Theorem \ref{thm:Twn^2=sumxi_ie_i} and Proposition \ref{thm:bauble}.
\end{proof}
}

Lastly, we rewrite Corollary \ref{cor:sq.wn.gamma} over $R$.

Recall the transformations $\T_s=q^{-1/2}T_s$ and $\xi=q^{1/2}-q^{-1/2}$.
So $\T_{w_n}^2=(q^{-\ell/2}T_{w_n})^2=q^{-\ell}T_{w_n}^2$.
We also have $\tilde\Gamma_\lambda=q^{-l_\lambda/2}\Gamma_\lambda$, and so $\xi^{l_\lambda}\tilde\Gamma_\lambda=q^{-l_\lambda/2}
(q^{1/2}-q^{-1/2})^{l_\lambda}\Gamma_\lambda=q^{-l_\lambda}(q-1)^{l_\lambda}\Gamma_\lambda$.   Recalling {from Proposition \ref{thm:bauble} that} $e_{n,i}=\sum_{l_\lambda=i}\Gamma_\lambda$,
the statement of Corollary \ref{cor:sq.wn.gamma} becomes
\begin{cor}\label{cor:sq.of.longest.word.q}
\[T_{w_n}^2=q^{\ell}\sum_{\lambda\vdash n}\left(1-q^{-1}\right)^{l_\lambda}\Gamma_\lambda.\]
\end{cor}

\section{The square root of the centre {of $\H$} }\label{sec:squareroot}

When $n=2$ the algebra $\H_n$ is commutative, so we restrict our attention in this section to $n\ge 3$.

\begin{defn}
  The \emph{square root} of $Z$, {denoted $\sqrt{Z}$},
  is defined {to be:}
  \[\sqrt{Z}:=\{h\in\H\mid h^2\in Z\}.\]
\end{defn}
Clearly $Z\subseteq\sqrt{Z}$.  {Also, $Z\neq\sqrt{Z}$ since} $T_{w_n}\in\sqrt{Z}\setminus Z$.  The fact that $T_{w_n}\in\sqrt{Z}\setminus Z$ is well-known and follows from the centrality of the square of the corresponding longest element in the braid group (of which $\H$ is a quotient), where that square is the full twist.  In this section we {define two additional elements which we show are in $\sqrt{Z}$,} and together with {$T_{w_n}$,} generate a commutative subalgebra of $\H$.
{In Theorem \ref{thm:Ri^2=elsymfn} we give} the forms of {the squares of these elements.}

We say an element $v\in\H$ is a {right} \emph{eigenvector} for multiplication by $h$ if $hv=kv$ for some {$k\in R$}.  {Left} eigenvectors are defined similarly.

{Recall the following well-known elements of {Z:
\begin{align*}
  x&:=\sum_{w\in S_n}T_w\\
  y&:=\sum_{w\in S_n}(-q)^{\ell-l(w)}T_w.
\end{align*}}

 These elements were defined in
{\cite[Section 3]{DipperJames86} and used there to study the permutation modules of $\H$.  They are widely used in the representation theory of {$\H$} (for example in \cite{DJ87}, \cite{Murphy92}, \cite{Murphy95}, and \cite{ErdmannNakano02} among many other places).  The modules $x\H$ and $y\H$ are called the trivial and alternating $\H$-modules respectively.} {Analogous elements can be defined for parabolic subalgebras of $\H$.}
{Both $x$ and $y$ are eigenvectors for the action of the generators of $\H$ (Lemma \ref{lem:props.of.x,y} \eqref{lem:sums.in.Z}), a property which they share with some other elements of $\sqrt Z$ (see Section \ref{sec:other.sq.roots}).}

The following result is well-known.

\begin{lem}{\ }

\begin{enumerate}\label{lem:props.of.x,y}
    \item\label{lem:sums.in.Z} For any $s\in S$, $T_sx=xT_s=qx$ and $T_sy=yT_s=-y$.  Moreover, $x$ and $y$ are central.
    \item  $xy=0$.
    \item\label{cor:gamma.forms.of.xn,yn} $x=\sum_{\lambda\vdash n}\Gamma_\lambda$ and $y=\sum_{\lambda\vdash n}(-q)^{\ell-l_\lambda}\Gamma_\lambda$.
    \item\label{cor:sq.of.xn,yn} $x^2=p(q)x$ and $y^2=(-1)^{\ell}p(q)y$.
\end{enumerate}
\end{lem}

We now define the elements $\bar x$ and $\bar y$, and investigate their properties in the remainder of this paper:
 \begin{align*}
  \bar x&:=x-T_{w_n}\\
  \bar y&:=y-T_{w_n}.
\end{align*}

\begin{lem}\label{lem:cross.terms} The following hold:
\begin{enumerate}
\item $\bar xT_{w_n}=T_{w_n}\bar x$,\label{eq:r1r3=r3r1}
\item $\bar yT_{w_n}=T_{w_n}\bar y$,\label{eq:r2r3=r3r2}
{ \item $\bar xT_{w_n}, \bar yT_{w_n}, \bar x\bar y\in Z$,\label{centrality} }
\item $\bar x\bar y=\bar y\bar x$.\label{eq:r1r2=r2r1}
\end{enumerate}
\end{lem}

\begin{proof}
The centrality of $x$ and $y$ from Lemma \ref{lem:props.of.x,y}\eqref{lem:sums.in.Z} implies \eqref{eq:r1r3=r3r1} and \eqref{eq:r2r3=r3r2}.  {Now ${w_n}$} acts on the set $S$ of generators {of $S_n$} via the graph automorphism $\rho_n$ defined {in Section \ref{sec:Prelim}}, {and this action extends to the Hecke algebra level as follows (remembering that $l(sw_n)=l(w_ns)=l(w_n)-1$ for any $s\in S$):
$T_{w_n}T_s=T_{\rho_n(s)}T_{\rho_n(s)w_n}T_s=T_{\rho_n(s)}T_{w_ns}T_s=T_{\rho_n(s)}T_{w_n}$.}
{Since $T_sx=xT_{s'}$ for any $s,s'\in S$, we have in particular $T_sx=xT_{\rho_n(s)}$. Hence, for any $s\in S$ (using the fact that $\rho_n$ is an involution), we have
\begin{align*}
\bar x T_{w_n}T_s   &=\bar x T_{\rho_n(s)}T_{w_n}\\
        &=\left(x-T_{w_n}\right) T_{\rho_n(s)}T_{w_n}\\
        &=\left(xT_{\rho_n(s)}-T_{w_n}T_{\rho_n(s)}\right) T_{w_n}\\
        &=\left(T_sx-T_sT_{w_n}\right) T_{w_n}\\
        &=T_s\left(x-T_{w_n}\right) T_{w_n}\\
        &=T_s\bar x T_{w_n}.
\end{align*}
Therefore, $\bar x T_{w_n} \in Z$. Similar arguments show that $\bar y T_{w_n}$ and $\bar x\bar y$ are central, proving (\ref{centrality}).}
Finally{,} the centrality of $x$ implies $\bar yx=x\bar y${, which} together with \eqref{eq:r2r3=r3r2} {implies} \eqref{eq:r1r2=r2r1}, and the proof is complete.
\end{proof}

\begin{prop}\label{prop:square.roots} {The elements $\bar x$, $\bar y$, $T_{w_n}$ are in $\sqrt{Z}\setminus Z$}.
\end{prop}
\begin{proof}
  As already mentioned, it is well known that $T_{w_n}\in\sqrt{Z}\setminus Z$. {By Lemma \ref{lem:props.of.x,y}\eqref{lem:sums.in.Z}, we have that $x,y\in Z$, from which we deduce that $\bar x,\bar y\notin Z$. Since $x^2$, $y^2$, $T_{w_n}^2 \in Z$, and since $T_{w_n}$ commutes with $\bar x$ and $\bar y$ (Lemma \ref{lem:cross.terms}), it follows that $\bar x^2, \bar y^2 \in Z$.}
\end{proof}

\begin{rem}
 {Note that while $\bar x+\bar y\notin Z$, some linear combinations of $\bar x$ and $\bar y$ are in $Z$ (for instance $\bar x-\bar y$), and all linear combinations of $\bar x$ and $\bar y$ are in $\sqrt{Z}$ since $\bar x^2, \bar y^2, \bar x\bar y\in Z$.}
\end{rem}
{The following is an immediate consequence of Lemma \ref{lem:cross.terms} and Proposition \ref{prop:square.roots}.
\begin{cor}\label{sqrtsubmodule}
 The $R$-submodule of $\H$ spanned by $\{\bar x,\bar y,T_{w_n}\}$ is contained in $\sqrt{Z}\setminus Z.$
\end{cor}}

\begin{prop}
 The elements $\bar x$ and $\bar y$ are not zero divisors.
\end{prop}
\begin{proof}
Suppose that $\bar xh=\left(x-T_{w_n}\right)h=0$ for some $h\in\H$.  Multiplying on the left by $x$ gives $\left(x^2-xT_{w_n}\right)h=0$, and therefore, using Lemma \ref{lem:props.of.x,y} parts \ref{lem:sums.in.Z} and \ref{cor:sq.of.xn,yn}, we have $\left(p(q)-q^{l(w_n)}\right)xh=0$.  This gives $xh=0$ since $\H$ is a free $R$-module, and subtracting $xh=0$ from $\left(x-T_{w_n}\right)h=0$ we obtain $T_{w_n}h=0$.
Thus $h=0$ since $T_w$ is invertible for all $w\in S_n$.
The proof for $\bar y$ is similar.
\end{proof}

\begin{prop}
{Let $\X$ be the $R$-subalgebra of $\H$ generated by {$\{\bar x,\bar y,T_{w_n}\}$}.
  Then $\X$ is commutative.}
\end{prop}
\begin{proof}
This is immediate from Lemma \ref{lem:cross.terms}.
\end{proof}

\begin{rem}
{Since $\bar x\bar yT_{w_n}\notin Z$, we have that $\bar x+\bar yT_{w_n}\notin\sqrt Z$. Hence, $\X\not\subseteq\sqrt{Z}$. This also shows that $\sqrt Z$ itself is not a subalgebra of $\H$. However, $\{\bar x^a\bar y^bT_{w_n}^c\mid a,b,c\in\N\}\subseteq \sqrt Z$, and, in fact, all elements of $\X$ with even exponent sum are central.}

{Unfortunately, $\sqrt{Z}$ (which includes $Z$) is not even
an $R$-module. Although we have shown that all $R$-linear combinations of the
particular elements $\bar x$, $\bar y$ and $T_{w_n}$ are contained in $\sqrt{Z}$ (Corollary \ref{sqrtsubmodule}), it is not true in general that other elements of $\sqrt{Z}$ have this property (see Section \ref{sec:other.sq.roots}),
and certainly unlikely that the sum of a non-central element in $\sqrt{Z}$ and an element of $Z$ will still be in $\sqrt{Z}$.}
\end{rem}

\begin{rem}
{There are no subalgebras of $\H$ which are contained in $\sqrt{Z}$ that contain elements in $\sqrt{Z}\setminus Z$.} This is because, {for instance,} if {$h\in\sqrt Z\setminus Z$, then $h+h^2\not\in\sqrt Z$.}
\end{rem}

\begin{thm}\label{thm:Ri^2=elsymfn}
  The squares of $\bar x$ and $\bar y$ can be expressed as linear combinations of elementary symmetric functions of Murphy elements {together with the identity}.  Specifically,
\begin{align*}
\bar x^2
&=\sum_{\lambda\vdash n}\left(p(q)-2q^{\ell}+q^{\ell-l_\lambda}\left(q-1\right)^{l_\lambda}\right)\Gamma_\lambda\\
&=\sum_{i=0}^{n-1}\left(p(q)-2q^{\ell}+q^{\ell-i}\left(q-1\right)^{i}\right) {e_{n,i}} ,\text{ and}\\
\bar y^2
&=\sum_{\lambda\vdash n}\left(-1\right)^{l_\lambda}q^{\ell-l_\lambda}\left(p(q)-2+(1-q)^{l_\lambda}\right)\Gamma_\lambda\\
&=\sum_{i=0}^{n-1}\left(-1\right)^{i}q^{\ell-i}\left(p(q)-2+(1-q)^i\right) {e_{n,i}}.
\end{align*}
\end{thm}
\begin{proof}
Immediately from the definition of $\bar x$ and $\bar y$ we have that
\begin{align*}
  \bar x^2&=x^2-2xT_{w_n}+T_{w_n}^2,\ \text{and}\\
  \bar y^2&=y^2-2yT_{w_n}+T_{w_n}^2.
\end{align*}
The forms for $x^2$, $y^2$ and $T_{w_n}^2$ from Lemma \ref{lem:props.of.x,y} parts \ref{cor:gamma.forms.of.xn,yn} and \ref{cor:sq.of.xn,yn} and Corollary \ref{cor:sq.of.longest.word.q}, together with Lemma \ref{lem:props.of.x,y} part \ref{lem:sums.in.Z}, then give that
\begin{align*}
\bar x^2&=p(q)x-2q^{\ell}x+ T_{w_n}^2\\
        &=\left(p(q)-2q^{\ell}\right)\sum_{\lambda\vdash n}\Gamma_\lambda+q^{\ell}\sum_{\lambda\vdash n}\left(1-q^{-1}\right)^{l_\lambda}\Gamma_\lambda.
\end{align*}
The expression for $\bar x^2$ in terms of $\{\Gamma_\lambda\mid\lambda\vdash n\}$ follows.  The corresponding expression for $\bar y^2$ is similarly derived. The expressions for both $\bar x^2$ and $\bar y^2$ in terms of $\{ {e_{n,i}} \mid 0\le i\le n-1\}$ follow immediately from {Proposition  \ref{thm:bauble}.}
\end{proof}

\section{Examples: {$\sqrt{Z}$ for $\H_3$ and $\H_4$} }\label{sec:other.sq.roots}

We have seen in Section \ref{sec:squareroot} that the commutative subalgebra $\X$ of $\H$ generated by $\{\bar x,\bar y,T_{w_n}\}$ is not contained in {$\sqrt{Z}$}.  {We see in the following examples that it is also the case that $\sqrt{Z}$ } is not contained in $\X$.

\subsection{ $\H_3$.}

{When $n=3$ it} is easy to check that while the elements $T_{s_1}-T_{s_2}$ and $T_{s_1s_2}-T_{s_2s_1}$ are in {$\sqrt{Z}$}, their span is not in {$\sqrt{Z}$} --- unlike $\bar x_3$, $\bar y_3$ and $T_{w_3}=T_{s_1s_2s_1}$.  However they do have the interesting property that they are eigenvectors for the multiplicative action of {$Z$.}  This parallels a similar fact for $x$ and $y$ (Lemma \ref{lem:props.of.x,y}\eqref{lem:sums.in.Z}). Recall from Section \ref{sec:Prelim} that $\{\G_{(1,1,1)}=T_1,\G_{(2,1)}=T_{s_1}+T_{s_2}+q^{-1}T_{s_1s_2s_1},\G_{(3)}=T_{s_1s_2}+T_{s_2s_1}+q^{-1}(q-1)T_{s_1s_2s_1}\}$
is an $R$-basis for {$Z$}. Then we have
\[\begin{array}{r@{}l@{\hspace{1mm}}c@{\hspace{1mm}}r@{}l}
    \G_{(2,1)}&(T_1-T_2)&=&(q-1)&(T_1-T_2)\\[1mm]
    \G_{(3)}&(T_1-T_2) &=&-q&(T_1-T_2)\\[1mm]
    \G_{(2,1)}&(T_{12}-T_{21})&=&(q-1)&(T_{12}-T_{21})\\[1mm]
    \G_{(3)}&(T_{12}-T_{21})&=&-q&(T_{12}-T_{21}).
\end{array}
\]

\begin{prop}
  When $n=3$, every element in $\sqrt{Z}\setminus Z$ can be written as a {$\Q(q)$--linear} combination of the following elements:
  \begin{align*}
    \bar x_3&:=T_1+T_{s_1}+T_{s_2}+T_{s_1s_2}+T_{s_2s_1}\\
    \bar y_3&:=T_1-q^{-1}(T_{s_1}+T_{s_2})+q^{-2}(T_{s_1s_2}+T_{s_2s_1})\\
    T_{w_3}&:=T_{s_1s_2s_1}\\
    R_4&:=T_{s_1}-T_{s_2}\\
    R_5&:=T_{s_1s_2}-T_{s_2s_1}.
  \end{align*}
\end{prop}
\begin{proof}
  Solving the system of equations created by evaluating \begin{multline*}(a_1T_1+a_2T_{s_1}+a_3T_{s_2}+a_4T_{s_1s_2}+a_5T_{s_2s_1}+a_6T_{s_1s_2s_1})^2T_s\\
=T_s(a_1T_1+a_2T_{s_1}+a_3T_{s_2}+a_4T_{s_1s_2}+a_5T_{s_2s_1}+a_6T_{s_1s_2s_1})^2\end{multline*}
for each $s\in S { =\{s_1,s_2\} } $ yields the following two solutions:
\begin{equation} \{a_2=a_3, \ a_4=a_5,\ a_6=q^{-1}a_3+q^{-1}(q-1)a_5\} \label{eq:A2.sol.central}\end{equation}
and
\begin{equation} \{a_1=-\frac{q-1}{2}(a_2+a_3)+\frac{q}{2}(a_4+a_5)\}.\label{eq:A2.sol.noncentral}
\end{equation}
The first solution, \eqref{eq:A2.sol.central}, characterizes the elements of {$Z$}: the relations given are exactly those given by \cite[Lemma (3.1)]{Fmb}. Thus the elements of {$\sqrt{Z}\setminus Z$} must satisfy equations \eqref{eq:A2.sol.noncentral}. The elements $\{\bar x_3,\bar y_3,T_{w_3},R_4,R_5\}$ {all satisfy these equations, and are easily seen} to be linearly independent. Thus they span a five-dimensional subspace of the solution space, which itself has five dimensions (there are five free parameters). {Hence} they span the solution space of \eqref{eq:A2.sol.noncentral}.
\end{proof}

\begin{cor}
  {For any subset $W$ of $\H$, let $\langle W \rangle$ denote the $R$-span of $W$. Then $\langle \{\bar x_3,\bar y_3,T_{w_3},R_4,R_5\} \rangle \cap Z_3=\langle \{x_3,y_3\} \rangle $.}
\end{cor}
\begin{proof}
  This is immediate from common solutions to \eqref{eq:A2.sol.central} and \eqref{eq:A2.sol.noncentral}.
\end{proof}

The squares of these five square roots of central elements (using Theorems \ref{thm:Ri^2=elsymfn} and \ref{thm:Twn^2=sumxi_ie_i}) are:
{\begin{align*}
  \bar x_3^2&=(2q^2+2q+1)\Gamma_{(1,1,1)}+(q+1)^2\Gamma_{(2,1)}+(3q+1)\Gamma_{(3)},\\
\bar y_3^2&=q^4(q^{2}+2q+2)\Gamma_{(1,1,1)}-q^3(q+1)^2\Gamma_{(2,1)}+q^3(q+3)\Gamma_{(3)},\\
  T_{w_3}^2 &=\Gamma_{(1,1,1)}+(1-q^{-1})\Gamma_{(2,1)}+(1-q^{-1})^2\Gamma_{(3)}\\
  R_4^2     &=2q\Gamma_{(1,1,1)}+(q-1)\Gamma_{(2,1)}-\Gamma_{(3)}\\
  R_5^2     &=-2q^2\Gamma_{(1,1,1)}-q(q-1){\Gamma_{(2,1)}} +q\Gamma_{(3)}.
\end{align*}
}

{Interestingly, we have the relation $R_5^2=-qR_4^2$. It follows that the central element given by $R_4^2$ has at least four distinct square roots over $\C(q^{1/2})$.}

\subsection{$\H_4$.}

As with $\H_3$, the algebra $\X$ generated by $\{\bar x_4,\bar y_4,T_{w_4}\}$ does not contain all elements of $\sqrt{Z}\setminus Z$.  Motivated by the relationships among the elements $\bar x_4$, $\bar y_4$ and the eigenvectors of $Z$ {when $n=4$}, as well as the experience in $\H_3$, we have identified several more elements of {$\sqrt{Z}$} by looking for eigenvectors of the multiplicative action of class elements.  The following are additional examples of square roots, which are clearly non-central since the coefficients of the shortest elements of conjugacy classes are not equal:

\begin{align*}
R_4&:=q\left(T_{s_1s_2}-T_{s_2s_1}+T_{s_3s_2}-T_{s_2s_3}\right)+(q-1)(T_{s_1s_3s_2}-T_{s_2s_1s_3})\\
    &\quad +T_{s_1s_2s_1s_3}-T_{s_1s_2s_3s_2}+T_{s_2s_3s_2s_1}-T_{s_1s_3s_2s_1}.\\
R_5&:=q^2(T_{s_1}+T_{s_3})+q(q-1)(T_{s_2s_1}+T_{s_2s_3}+T_{s_1s_3})+(q-1)^2T_{s_2s_1s_3}\\
    &\quad -q(T_{s_1s_2s_1}+T_{s_2s_3s_2}+T_{s_1s_2s_3}+T_{s_3s_2s_1})-(q-1)(T_{s_1s_2s_1s_3}\\
    &\quad +T_{s_2s_3s_2s_1}+T_{s_2s_1s_3s_2})+T_{s_1s_2s_1s_3s_2}+T_{s_2s_1s_3s_2s_1}.\\
R_6&:=q^2T_{s_2}+q(q-1)(T_{s_1s_2}+T_{s_3s_2})-q(T_{s_1s_2s_1}+T_{s_2s_3s_2}+T_{s_1s_2s_3}\\
    &\quad +T_{s_3s_2s_1}-T_{s_2s_1s_3})+(q^2-q+1)T_{s_1s_3s_2}-(q-1)(T_{s_1s_2s_3s_2}\\
    &\quad +T_{s_1s_3s_2s_1})+T_{s_1s_2s_3s_2s_1}.\\
\end{align*}

The set $\{\bar x_4,\bar y_4, T_{w_4}, R_4, R_5, R_6\}$ is linearly independent.  While the elements $\{\bar x_4,\bar y_4, T_{w_4}\}$ commute with each other and with elements of $\{R_4, R_5, R_6\}$, the elements of $\{R_4, R_5, R_6\}$ do not commute among themselves.

Unlike the $\H_3$ case, we cannot use a dimension argument here to prove that we have all square roots.  In fact{,} as the dimension of the algebra here is 24{,} it is rather unlikely that this set spans all square roots.

\bibliographystyle{line}

\end{document}